\documentclass[11pt]{article}
\usepackage[a4paper, total={5in, 9in}]{geometry}

\usepackage{amsmath,amssymb,dsfont,amsfonts,amsthm,amscd}
\usepackage{tabularx}
\usepackage{graphicx}
\usepackage{times}
\usepackage{xcolor}
\usepackage{float}
\usepackage[pdfencoding=auto]{hyperref}
\usepackage{url}
\usepackage{array}
\usepackage{bm}
\usepackage{authblk}
\usepackage[numbers]{natbib}
\usepackage{setspace}
\usepackage{booktabs}
\usepackage{stackengine}
\usepackage{ulem}
\usepackage{rotating}
\usepackage[caption=false]{subfig}
\stackMath

\usepackage{tcolorbox}
\usepackage{dsfont}
\usepackage{booktabs}
\usepackage{xcolor}
\definecolor{dgreen}{rgb}{0.0, 0.5, 0.0}
\usepackage{float}
\usepackage{algorithm}
\let\classAND\AND
\let\AND\relax
\usepackage{algorithmic}

\let\AND\classAND
\AtBeginEnvironment{algorithmic}{\let\AND\algoAND}

\title{Optimal Risk Scores for Continuous Predictors}

\author[]{Cristina Molero-R\'io \thanks{Corresponding author: molero@lix.polytechnique.fr, mmolero@us.es}\,\,}
\author[]{Claudia D'Ambrosio}
\affil[]{LIX CNRS, École Polytechnique, Institute Polytechnique de Paris,\\ Route de Saclay, Palaiseau, 91128, France}

\date{}

\providecommand{\keywords}[1]{\textbf{\textit{Keywords---}} #1}

\usepackage{pdfpages}

\begin{document}



\maketitle

\begin{abstract}

In this paper, we propose a novel Mixed-Integer Non-Linear Optimization formulation to construct a risk score, where we optimize the logistic loss with sparsity constraints. Previous approaches are typically designed to handle binary datasets, where continuous predictor variables are discretized in a preprocessing step by using arbitrary thresholds, such as quantiles. In contrast, we allow the model to decide for each continuous predictor variable the particular threshold that is critical for prediction. 
The usefulness of the resulting optimization problem is tested in synthetic datasets.
\end{abstract}


\section{Introduction}
\label{Introduction}

With the widespread use of Machine Learning applications that guide decisions in many daily tasks, interpretability has become essential for trust in them \cite{rudin2022}. This is supported by the so-called right-to-explanation in algorithmic decision-making, imposed by the European Union since 2018 \cite{goodman2017}, that encourages to design models that are understandable and enable explanation to the users. Mathematical Optimization has been shown to be a flexible enough tool for dealing with this necessity and creating Interpretable Machine Learning models, see \cite{gambella2021} for a survey. Some examples include the construction of decision trees \cite{carrizosa2021mathematical}, rule sets \cite{lawless2023interpretable}, and scoring systems \cite{ustun2016supersparse}.

Another fundamental tool for Interpretable Machine Learning models is Risk Scores \cite{rudin2022}. They are widely used to support human decision-making in relevant domains such as Criminal Justice \cite{latessa2010creation} to assess the risk of recidivism when sentencing, Finance \cite{siddiqi2017intelligent} to assess the risk of default on a loan, or Healthcare \cite{moreno2005saps} to assess the risk of mortality of ICU patients, among others. Risk Scores consist of logistic regression models applied to scoring systems.

A logistic regression model is a statistical model to predict a binary prediction variable of the form $Y\in\left\lbrace 0,1\right\rbrace$ that assumes that the log-odds are a linear function of a set of independent variables $X=\left(X_1,\ldots, X_p\right)\in\mathbb{R}^p$, that is, 
\begin{equation}
\label{eq:logit_model}
\log\dfrac{P\left(Y=1\,|\,X\right)}{P\left(Y=0\,|\, X\right)} =  w_0 + \sum_{j=1}^p w_j X_j = w_0 + X^\top\bm{w},
\end{equation} for some $w_0\in\mathbb{R}, \,\, \bm{w}= \left(w_1,\ldots,w_p\right)\in\mathbb{R}^p$, where $ P\left(Y=0\,|\,X\right) = 1 - P\left(Y=1\,|\, X\right)$. Inverting Equation \eqref{eq:logit_model}, one can derive the expressions of the membership probability to each category as follows:
\begin{equation}
\nonumber
  \dfrac{P\left(Y=1\,|\,X\right)}{1-P\left(Y=1\,|\,X\right)} =   \exp\left( w_0 + X^\top\bm{w}\right)\\
\end{equation}
with
\begin{align}
\begin{split}
  \label{eq:prob1}    
P\left(Y=1\,|\,X\right) &= \dfrac{\exp\left( w_0 + X^\top\bm{w}\right)}{1+\exp\left( w_0 + X^\top\bm{w}\right)} = 
 \dfrac{1}{1+\exp\left( -\left( w_0 + X^\top\bm{w}\right)\right)}\\
 & = \sigma\left(w_0 + X^T\bm{w}\right)
 \end{split}
\end{align} and
\begin{align}
\begin{split}
\label{eq:prob2}
    P\left(Y=0\,|\,X\right) & = \dfrac{\exp\left( -\left( w_0 + X^\top\bm{w}\right)\right)}{1+\exp\left( -\left( w_0 + X^\top\bm{w}\right)\right)} = \dfrac{1}{1+\exp\left(  w_0 + X^\top\bm{w}\right)} \\
    &=\sigma\left(  -\left( w_0 + X^\top\bm{w}\right)\right),
\end{split}
\end{align} where $\sigma: \mathbb{R} \longrightarrow\left[ 0,1\right]$ is the sigmoid function.

A scoring system is a linear classifier with integer coefficients. Scoring systems are commonly based on binary or indicator predictors in such a way that the user is required to perform simple and few calculations with integer numbers in order to compute a score and make a prediction.

\begin{figure}[!ht]
    \begin{center}
        \includegraphics[scale=0.3]{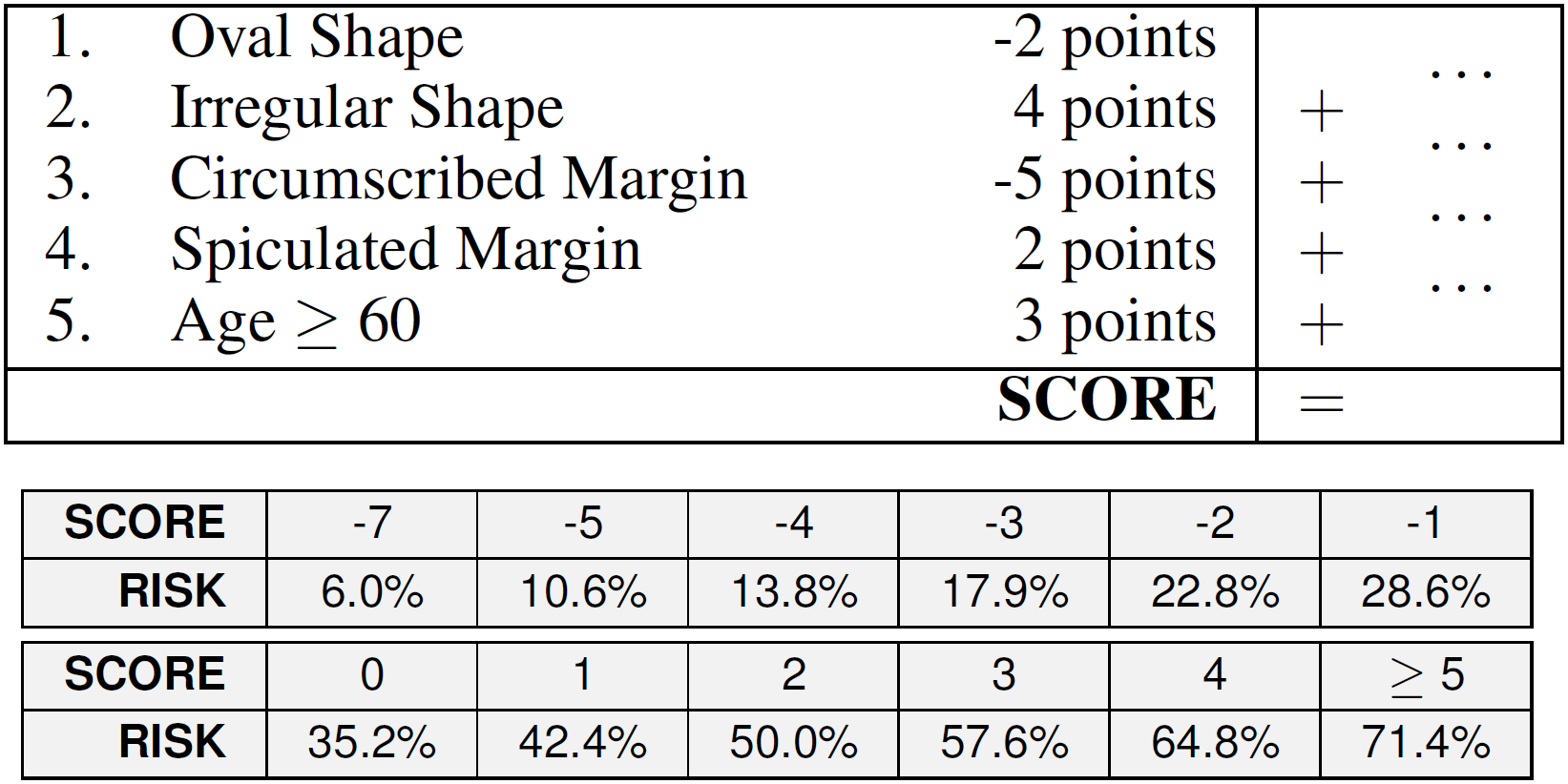}
        \caption[Fig]{Risk score from \cite{liu2022} on the \texttt{mammo} dataset, which consists of a sample of biopsy patients. The model tries to predict the risk of malignancy of a breast lesion.}
        \label{logoRoadef}
    \end{center}
\end{figure}

Consequently, after applying logistic regression to the scoring system, risk scores will produce a discrete set of probabilities, each of them based on the possible integer outputs coming from the scoring system, as illustrated in Figure \ref{logoRoadef}. It shows an example of a risk score for an application in healthcare, corresponding to the well-known \texttt{mammo} \cite{elter2007prediction} dataset. The \texttt{mammo} dataset consists of 14 characteristics measured on 961 biopsy patients. This particular dataset only includes binary features, although this is not the case for others in the domain such as \texttt{breastcancer} \cite{mangasarian1995breast} that contain continuous ones as well. The goal is to predict the severity – benign or malignant – of a mammographic mass lession, according to such characteristics. The box at the top represents the scoring system, which is $[\textup{Scoring System}] := \textup{SCORE} \,\,+\,\, w_0$ where $\textup{SCORE} := -2\cdot [\textup{Oval Shape}]\,\, +\,\, 4 \cdot [\textup{Irregular Shape} ] \,\,-\,\,5 \cdot [\textup{Circumscribed Margin}]  \,\,+\,\, 2 \cdot [\textup{Spiculated Margin}] \,\,+\,\, 3 \cdot [\textup{Age}\geq60 ]$. Note that the five features are indicator, taking value 1 when what is written in [] is true, and $w_0$ is the independent term. The table at the bottom is the conversion of the score to the probability of being assigned to the positive class, in this case, the malignant category. In other words, $[\textup{Risk Score}] :=\sigma\left(\dfrac{1}{m}[\textup{Scoring System}]\right) =\sigma\left(\dfrac{1}{m}[\textup{SCORE} \,\,+\,\, w_0]\right) = \textup{RISK}$, for some scaling parameter $m$. As a  result, the user can easily determine the risk of malignancy of a breast lesion by adding points for its shape, its margin, and the patient’s age. If the score is above a threshold, the patient would be recommended to further tests. Suppose a 65 years old patient presents a mammographic mass lession with oval - and not irregular - shape and with spiculated - and not circumscribed - margin. Then, their SCORE is $-2\cdot 1 + 4 \cdot 0 -5 \cdot 0 + 2 \cdot 1 + 3 \cdot 1 = 3$, which is translated into a RISK of malignant severity of $57.6\%$.

The problem of constructing risk scores dates back to the 1930's, see \cite{burgess1928factors}.  Risk scores have traditionally been constructed either by domain experts or heuristically \cite{bmj,antman2000timi}, achieving a low computational cost, at the expense of prediction performance. Recent work has aimed to fully optimize risk scores \cite{ustun2019,liu2022}. While efficient, these approaches are designed to work with binary predictors. Thus, continuous predictors are discretized in a preprocessing step by using arbitrary thresholds, such as quantiles. Our motivation is to decide the optimal threshold to the corresponding continuous predictor while constructing the model, and, thus, hopefully having a more tailored binarized continuous predictor. For that, we propose a Mixed Integer Non Linear Optimization (MINLO) formulation that seeks a good trade-off between prediction performance and interpretability.

\section{Related work}
\label{Related work}

In this section, we review the different approaches that have been considered in the past to solve the problem of learning risk scores.

Traditionally, risk scores have been built by a panel of experts in the domain \cite{bmj,gage2001validation}, who designed the model with no data at hand and, thus, no notion on its predictive performance, potentially yielding to suboptimal models \cite{ustun2019}.

Since the integrality constraint on the scoring coefficients hinders the tractability of the problem, there has been a part of the literature devoted to implement select-regress-and-round schemes \cite{jung2017simple}. These approaches aim at fitting real scoring coefficients using the continuous relaxation of the problem with a preliminary feature selection or some kind of regularization term, followed by a rounding heuristic \cite{chevaleyre2013rounding,cole1993algorithm} in order to obtain integer-valued scoring coefficients. The traditional rounding to the nearest integer is used in \cite{cole1993algorithm}. A similar approach is developed in this paper in order to obtain a first feasible solution of the problem. Rounding can sometimes go against the loss gradient and harm predictive performance. There have been some attempts to improve the traditional rounding. In \cite{chevaleyre2013rounding}, they investigate two rounding schemes: a randomized one, where the scoring coefficients are rounded up or down randomly, and a sequential one, where 
the randomized scheme is applied sequentially. In \cite{ustun2019}, they also implement a sequential rounding, making the best choice, up or down, according to the regularized
logistic loss.

More recently, several optimization-based approaches \cite{ustun2017optimized,ustun2019,liu2022} have been proposed to fully optimize risk scores. In \cite{ustun2017optimized,ustun2019}, a mixed-integer non-linear program that minimizes the logistic loss, penalizes the $\ell_0$-norm for sparsity, and restricts coefficients to small integers, is proposed. A cutting plane algorithm to solve it is presented, where a surrogate problem with a linear approximation of the loss function is iteratively solved. The method shows to scale linearly with the number of samples in the dataset. In \cite{liu2022}, a more efficient algorithm to solve the problem in  \cite{ustun2017optimized,ustun2019} is proposed, where the $\ell_0$-penalization term is translated into a budget constraint. We use this formulation as starting point of our approach, and extend it to handle continuous predictors. Their algorithm consists of three steps. First, the continuous relaxation of the problem is solved by using beam search for best subset selection. Second, a Rashomon set is constructed around this continuous solution, that is, a pool of almost equally accurate solutions is constructed by
swapping one feature at a time, and comparing the new loss with a prefixed tolerance gap level. In the third step, for each coefficient vector in the pool, a risk score with integer scoring coefficients is obtained by introducing multipliers and a customized rounding scheme.

Bayesian methods have been also proposed for building scoring systems that make use of a Markov Chain Monte Carlo procedure for obtaining their solution. Such is the case of \cite{ertekin2015bayesian}, where Bayesian priors control the shape of the coefficients, in contrast to regularization penalties.

Other losses, different from the logistic we use, such as the 0-1 loss \cite{ustun2016supersparse}, or the hinge loss \cite{billiet2018interval,sokolovska2017fused}, and other link functions, such as the probit \cite{prentice1976generalization}, have been proposed to address similar questions. Apart from optimizing prediction performance and sparsity, there are works that focus on fairness, such as \cite{petersen2023assessing}.

\section{Optimal risk scores}
\label{Optimal risk scores}

In this section, we introduce our approach step by step. Section \ref{Introduction Section 3} motivates the use of the logistic loss as objective function. Section \ref{Original formulation for binary input} introduces the state-of-the-art formulation for binary inputs that we depart from. Section \ref{Reformulation for binary input} proposes a reformulation that will be used in our approach. Section \ref{New formulation with bilinear terms} introduces a first sketch of our formulation for continuous inputs using bilinear terms. Finally, Section \ref{New formulation with linear terms} eliminates these bilinear terms and introduces our final proposal.

\subsection{Introduction}
\label{Introduction Section 3}

In logistic regression, the coefficients are commonly estimated via maximum likelihood, which is equivalent to minimizing the logistic loss, as we show below.

According to Equations \eqref{eq:prob1} and \eqref{eq:prob2}, the conditional distribution is given by 
\begin{equation}
    \label{eq:conditional_distribution}
    P\left(Y=y | X=\bm{x}\right) = \sigma\left(w_0 + \bm{x}^T\bm{w}\right)^y \sigma\left(-\left( w_0 + \bm{x}^T\bm{w}\right)\right)^{1-y}.
\end{equation}
Given $\left\lbrace \left(\bm{x}_i,y_i\right)\right\rbrace_{1\leq i\leq N}$, an i.i.d. training set from $X\times Y$ of size $N$, we want to maximize the likelihood:
\begin{align}
    \max\limits_{w_0,\bm{w}}\,\, & \prod_{i=1}^N P\left(Y=y_i | X=\bm{x}_i\right) =\\
 \max\limits_{w_0,\bm{w}}\,\, & \prod_{i=1}^N\sigma\left(w_0 + \bm{x}_i^\top\bm{w}\right)^{y_i} \sigma\left(-\left( w_0 + \bm{x}_i^\top\bm{w}\right)\right)^{1-y_i};
\end{align}
which is equivalent to maximizing the log-likelihood:
\begingroup
\begin{align}
    \max\limits_{w_0,\bm{w}}\,\, & \sum_{i=1}^N y_i \log \left(\sigma\left(w_0 + \bm{x}_i^\top\bm{w}\right)\right) + \left( 1-y_i\right)\log\left(\sigma\left(-\left( w_0 + \bm{x}_i^\top\bm{w}\right)\right)\right) =\\
     \max\limits_{w_0,\bm{w}}\,\, & \sum_{i: y_i = 1} \log \left(\sigma\left(w_0 + \bm{x}_i^\top\bm{w}\right)\right) + \sum_{i: y_i = 0} \log\left(\sigma\left(-\left( w_0 + \bm{x}_i^\top\bm{w}\right)\right)\right) =\\
     \max\limits_{w_0,\bm{w}}\,\, & \sum_{i: y_i = 1} -\log \left(1+\exp\left( -\left( w_0 + \bm{x}_i^\top\bm{w}\right)\right)\right) \notag \\
     & \hspace{4.5cm} + \sum_{i: y_i = 0} -\log \left(1+\exp\left(  w_0 + \bm{x}_i^\top\bm{w}\right)\right);
\end{align}
\endgroup which is equivalent to the minimization of the logistic loss for new class labels  $\tilde{y}_i$ where $\tilde{y}_i=1$ if $y_i=1$, and $\tilde{y}=-1$ if $y_i=0$, $\forall i = 1,\ldots, N$:
\begingroup
\begin{equation}
\label{logisticloss}
    \min\limits_{w_0,\bm{w}}\,\, \sum_{i = 1}^N \log\left(1+\exp\left(-\tilde{y}_i\left( w_0 + \bm{x}_i^\top\bm{w}\right)\right)\right).
\end{equation}
\endgroup
With abuse of notation, $y_i$ will replace $\tilde{y}_i$ in the following.

\subsection{Original formulation for binary input}
\label{Original formulation for binary input}

As said in Section \ref{Introduction}, we depart from the formulation in \cite{liu2022} that was proposed to construct a risk score with binary input, that is, $\bm{x}_i\in\left\lbrace 0,1\right\rbrace^p$, and prediction $y_i\in\left\lbrace -1,1\right\rbrace$. The formulation reads as follows:
\begin{align}
    \min\limits_{w_0,\bm{w}}\,\, & \sum\limits_{i=1}^N \log\left( 1+\exp\left(-y_i\left(\bm{x}_i^\top\bm{w} + w_0\right)\right)\right) \label{eq:original_obj}\\
    \textup{s.t.} \,\,& \lVert \bm{w} \rVert_0 \leq k \label{eq:original_sparsity}\\
    & w_j \in \left[-5,5\right], w_j\in\mathbb{Z}, j=1,\ldots,p  \label{eq:original_small}\\
    & w_0 \in \mathbb{Z} \label{eq:original_intercept}.
\end{align}
The objective function \eqref{eq:original_obj} minimizes the logistic loss. Constraint \eqref{eq:original_sparsity} induces sparsity by restricting the number of binary features to be used in the model to a predefined and small value for parameter $k$. Constraints \eqref{eq:original_small} impose integrality to the scoring coefficients, as well as box constraints that are user-defined and can be different for each coefficient, but, for simplicity, 5 is used. Constraint \eqref{eq:original_intercept} assumes the integrality of the intercept as well.

\subsection{Reformulation for binary input}
\label{Reformulation for binary input}

In order to deal with the $\ell_0$-norm appearing in Constraint \eqref{eq:original_sparsity}, we define new binary decision variables $\alpha_j,\,\, j=1,\ldots,p$, that will take value $1$ if feature $j$ is used in the model and 0, otherwise. Note that, the introduction of variables $\alpha$ allows us to write the sparsity constraint linearly, at the cost of introducing a binary variable per feature.
Thus:
\begin{align}
    \min\limits_{w_0,\bm{w},\bm{\alpha}}\,\, & 
    \eqref{eq:original_obj} \notag\\
    \textup{s.t.} \,\,& \eqref{eq:original_small}-\eqref{eq:original_intercept} \notag\\
    & \sum_{j=1}^p \alpha_j \leq k \label{eq:reformulation_alpha}\\
    & -5\alpha_j\leq w_j\leq 5\alpha_j,\,\, j=1,\ldots,p\label{eq:reformulation_alpha_w}\\
    & \alpha_j \in \left\lbrace 0,1\right\rbrace, j = 1,\ldots,p\label{eq:int_alpha}
\end{align}

Constraint \eqref{eq:reformulation_alpha} replaces \eqref{eq:original_sparsity} with the inclusion of the new binary decision variables defined by Constraints \eqref{eq:int_alpha}. Constraints \eqref{eq:reformulation_alpha_w} link $\alpha_j$ and $w_j$, by imposing $w_j$ to be $0$ whenever feature $j$ is not used in the model, that is, $\alpha_j=0$, and simple bounds, otherwise. The objective function \eqref{eq:original_obj} and constraints \eqref{eq:original_small}-\eqref{eq:original_intercept} remain the same.

\subsection{New formulation with bilinear terms}
\label{New formulation with bilinear terms}

We now introduce the first sketch of our contribution, i.e., an extension of the formulation which considers continuous inputs instead. Similarly, for the sake of interpretability, these continuous predictors are being binarized at the same time the model is constructed.

Let us assume that we depart from a sample $\left\lbrace \left(\bm{x}_i,y_i\right)\right\rbrace_{i=1}^N$ from $X\times Y$, where $X_j\in\left[\underline{x}_j,\overline{x}_j\right],\,\, j = 1,\ldots,p,$ and $Y\in\left\lbrace -1, 1\right\rbrace$.

First, we need to decide which conditions we consider to binarize the continuous predictors. We choose to implement conditions with the following shape: $X_j\leq t_j$ and $X_j\geq t_j$, where $t_j\in\left[\underline{x}_j,\overline{x}_j\right],\,\, j=1,\ldots,p,$ are new continuous decision variables that will represent the threshold values for predictors $j=1,\ldots,p$, respectively. One should not care about the case $X_j=t_j$; it would consist of just summing up two values in case they both are selected, and handling with no strict inequalities will ease the formulation.

We need to include the following new decision variables: 

$w^{\leq}_j \in [-5, 5],\,\, w^{\leq}_j\in\mathbb{Z},\,\,j=1,\ldots,p$: scoring coefficient corresponding to the condition $X_j\leq t_j$;

$w^{\geq}_j \in [-5, 5],\,\, w^{\geq}_j \in\mathbb{Z},\,\,j=1,\ldots,p$: scoring coefficient corresponding to the condition $X_j\geq t_j$.

In order to compute properly the scoring, we need to define new decision variables to specify which condition is satisfied by each observation, that is:

$b_{ij}\in\left\lbrace 0,1\right\rbrace,\,\, i=1,\ldots,N,\,\,j=1,\ldots,p$: binary decision variable that takes value $1$ if $x_{ij} \leq t_j$ and 0, if $x_{ij}\geq t_j$.

Note that variables $b$ are necessary to model the score formula, i.e., to make sure that a given score is added only if the feature is $\leq$ (or $\geq$) of the computed threshold $t$.




Then, the model will read as follows. 

\begingroup
\begin{align}
    \min\limits_{\substack{w_0,\\
    \bm{w^\leq}, \bm{w^\geq},\\
    \bm{\alpha},\bm{t},\bm{b}}}\,\, & \sum\limits_{i=1}^N \log\left( 1+\exp\left(-y_i\left(\sum_{j=1}^p \left( w^{\leq}_{j}b_{ij} + w^{\geq}_{j}\left(1-b_{ij}\right) \right) + w_0\right)\right)\right)\label{eq:bilinear_obj}\\
    \textup{s.t.} \,\,& \eqref{eq:original_intercept}-\eqref{eq:reformulation_alpha}, \eqref{eq:int_alpha}\notag\\
    &  -5\alpha_j\leq w^{\leq}_j\leq 5\alpha_j,\,\, j=1,\ldots,p \label{eq:bilinear_alpha_w_1}\\
    &  -5\alpha_j\leq w^{\geq}_j\leq 5\alpha_j,\,\, j=1,\ldots,p \label{eq:bilinear_alpha_w_2}\\
    & (\underline{x}_j - \overline{x}_j) b_{ij} + t_j \leq x_{ij} \leq t_j + (\overline{x}_j-\underline{x}_j) (1-b_{ij}), \,\, i = 1,\ldots,N, \label{eq:bilinear_thresholds}\\
    & \hspace{7.5cm} j=1,\ldots,p, \notag\\
    &  w^{\leq}_j,w^{\geq}_j\in\mathbb{Z}, j=1,\ldots,p\label{eq:bilinear_w}\\
    & t_j \in \left[\underline{x}_j,\overline{x}_j\right], j = 1,\ldots,p,\label{eq:bilinear_tj}\\
    & b_{ij} \in \left\lbrace 0,1\right\rbrace, i = 1,\ldots,N, j = 1,\ldots,p\label{eq:bilinear_bij}
\end{align}
\endgroup

In the objective function \eqref{eq:original_obj}, we were adding the weight $w_j$ every time $x_{ij}=1$. In this formulation, we have the binary decision variables $b_{ij}$ that play the role of the binarization, in such a way that the objective function \eqref{eq:bilinear_obj} will consider $w_j^\leq$ every time $b_{ij}=1$ and $w_j^\geq$, otherwise. Constraints \eqref{eq:original_intercept}-\eqref{eq:reformulation_alpha} and \eqref{eq:int_alpha} remain the same. Constraints \eqref{eq:bilinear_alpha_w_1} and \eqref{eq:bilinear_alpha_w_2} replace Constraints \eqref{eq:reformulation_alpha_w} and link, similarly, $\alpha_j$ with $w_j^\leq$ and $w_j^\geq$. Constraints \eqref{eq:bilinear_thresholds} ensure that the thresholds are well defined. If $b_{ij}=0$, then $x_{ij} \geq t_j$. If $b_{ij}=1$, then $x_{ij} \leq t_j$. Constraints \eqref{eq:bilinear_w}, \eqref{eq:bilinear_tj}, and \eqref{eq:bilinear_bij} define the range of the new decision variables.

\subsection{New formulation with linear terms}
\label{New formulation with linear terms}

In the following, we formulate our final proposal. We will introduce new decision variables in order to replace the bilinear terms appearing in the objective function in the previous formulation, that is:

$\beta_{ij}^{\leq} = w^{\leq}_{j}b_{ij} \in [-5, 5],\,\,i=1,\ldots,N,\,\, j=1,\ldots,p$: term added to the score in the objective function if $x_{ij}\leq t_j,\,\,i=1,\ldots,N,\,\, j=1,\ldots,p$;

$\beta_{ij}^{\geq} = w^{\geq}_{j}\left(1-b_{ij}\right)  \in [-5, 5],\,\,i=1,\ldots,N,\,\, j=1,\ldots,p$: term added in the objective function if $x_{ij}\geq t_j,\,\,i=1,\ldots,N,\,\, j=1,\ldots,p$;

\noindent with additional constraints $\forall i =1,\ldots,N,\,\, \forall j = 1,\ldots,p$:

\bigskip
\begin{tabular}{l l}
      $\beta_{ij}^\leq \leq 5b_{ij}$ & \hspace{1cm} $\beta_{ij}^\geq \leq 5\left(1-b_{ij}\right)$\\
     $\beta_{ij}^\leq \leq w_j^\leq - 5 b_{ij} + 5$  & \hspace{1cm} $\beta_{ij}^\geq \leq w_j^\geq - 5 (1-b_{ij}) + 5 = w_j^\geq + 5b_{ij}$\\
    $\beta_{ij}^\leq \geq - 5b_{ij}$ & \hspace{1cm} $\beta_{ij}^\geq \geq - 5\left(1-b_{ij}\right)$\\
    $\beta_{ij}^\leq \geq 5b_{ij} + w_j^\leq -5 $   & \hspace{1cm} $\beta_{ij}^\geq \geq 5\left(1-b_{ij}\right) + w_j^\geq -5 = -5b_{ij} + w_j^\geq$
\end{tabular}
\bigskip

\noindent which represent the Fortet's reformulation \cite{fortet1960applications} of the given bilinear functions. This reformulation is exact given that $b_{ij}$ are binary decision variables. In particular, when $b_{ij}=0$, the first and the third constraints on the left impose $\beta_{ij}^\leq = 0$, while the second and the forth are redundant. On the contrary, when $b_{ij}=1$, the second and the forth impose that $\beta_{ij}^\leq = w_j^\leq$, while the first and the third are redundant. The four constraints on the right work similarly but for the bilinear term $w_j^\geq (1-b_{ij})$. Note that, thanks to the introduction of additional variables $\beta_{ij}^\leq$ and $\beta_{ij}^\geq$, we managed to linearize the products $w^{\leq}_{j}b_{ij}$ and $w^{\geq}_{j}\left(1-b_{ij}\right)$, respectively. We decided to linearlize these products to keep the objective convex and, consequently, simplify the formulation at cost of additional variables and constraints.
Thus, we propose the following formulation, that we fully describe here:

\begingroup
\allowdisplaybreaks
\begin{align}
    \min\limits_{\substack{w_0,\\
    \bm{w^\leq}, \bm{w^\geq},\\ \bm{\alpha},\bm{t},\bm{b},\\ \bm{\beta^\leq},\bm{\beta^\geq}}}\,\, & \sum\limits_{i=1}^N \log\left( 1+\exp\left(-y_i\left(\sum_{j=1}^p \left( \beta^{\leq}_{ij} + \beta^{\geq}_{ij} \right) + w_0\right)\right)\right) \label{eq:linear_obj}\\
    \textup{s.t.} \,\,& \sum_{j=1}^p \alpha_j \leq k\\
    & -5\alpha_j\leq w^{\leq}_j\leq 5\alpha_j,\,\,  j=1,\ldots,p\\
    & -5\alpha_j\leq w^{\geq}_j\leq 5\alpha_j,\,\, j=1,\ldots,p\\
    & (\underline{x}_j - \overline{x}_j) b_{ij} + t_j \leq x_{ij} \leq t_j + (\overline{x}_j-\underline{x}_j) (1-b_{ij}), \,\,i = 1,\ldots,N,\\
    & \hspace{7.5cm} j=1,\ldots,p \notag\\
    & \beta_{ij}^\leq \leq 5b_{ij}, \,\,i = 1,\ldots,N, j=1,\ldots,p\\
    & \beta_{ij}^\leq \leq w_j^\leq - 5 b_{ij} + 5, \,\,i = 1,\ldots,N, j=1,\ldots,p\\
    & \beta_{ij}^\leq \geq - 5b_{ij}, \,\,i = 1,\ldots,N, j=1,\ldots,p\\
    & \beta_{ij}^\leq \geq 5b_{ij} + w_j^\leq -5, \,\,i = 1,\ldots,N, j=1,\ldots,p\\ 
    & \beta_{ij}^\geq \leq 5\left(1-b_{ij}\right), \,\,i = 1,\ldots,N, j=1,\ldots,p\\
    & \beta_{ij}^\geq \leq w_j^\geq + 5b_{ij}, \,\,i = 1,\ldots,N, j=1,\ldots,p\\
    & \beta_{ij}^\geq \geq - 5\left(1-b_{ij}\right), \,\,i = 1,\ldots,N, j=1,\ldots,p\\
    & \beta_{ij}^\geq \geq  -5b_{ij} + w_j^\geq, \,\,i = 1,\ldots,N, j=1,\ldots,p\\
    & w^{\leq}_j,w^{\geq}_j\in\mathbb{Z}, \,\,j=1,\ldots,p\\
    & w_0 \in \mathbb{Z},\\
    & \alpha_j \in \left\lbrace 0,1\right\rbrace, \,\, j=1,\ldots,p\\
    & t_j \in \left[\underline{x}_j,\overline{x}_j\right], \,\, j=1,\ldots,p\\
    & b_{ij} \in \left\lbrace 0,1\right\rbrace, \,\,i = 1,\ldots,N, j=1,\ldots,p\\
     & \beta^{\leq}_{ij},\beta^{\geq}_{ij}\in\mathbb{Z}, \,\,i = 1,\ldots,N, j=1,\ldots,p, \label{eq:linear_last}
\end{align}
\endgroup
which reads as a Mixed-Integer Non-Linear Optimization problem with linear constraints. The objective is non-linear while convex. Apart from box constraints, there are $1+(2+9N)p$ linear constraints, $p(1+N)$ binary, $1+2p(1+N)$ integer and $p$ continuous decision variables.

\section{A simple heuristic}
The formulations presented in the previous sections belong to the class of convex mixed integer nonlinear programs. These problems are difficult to solve for their combinatorial nature. However, we know efficient algorithms to solve their continuous relaxations. Thus, we devise a simple heuristic with the aim of finding a first feasible solution. We present it in Algorithm \ref{alg:heur}.
\begin{algorithm}[h!]
\begin{algorithmic}[1] 
\STATE Solve the continuous relaxation of formulation \eqref{eq:linear_obj}-\eqref{eq:linear_last} and let $\overline{\alpha}, \overline{\beta^{\leq}}, \overline{\beta^{\geq}}, \overline{b}, \overline{w_0}, \overline{w^\leq}, \overline{w^\geq}, \overline{t}$ be one of its optimal solutions.
\STATE $\overline{w_{0}} = \lfloor \overline{w_{0}} \rceil$ \label{algo1:round_w0}
\STATE Consider $\overline{\alpha}$ in non-increasing order and set the first $k$ to $1$. \label{algo1:round_alphas}
\FOR{$j = 1; j<=p; j++$}
    \IF{$\overline{\alpha}_j = 1$}
        \STATE $\overline{w_{j}^\leq} = \lfloor \overline{w_{j}^\leq} \rceil$
        \STATE $\overline{w_{j}^\geq} = \lfloor \overline{w_{j}^\geq} \rceil$
        \FOR{$i = 1; i\leq N; i++$}
            \IF{$x_{ij} \leq \overline{t_j}$}
                \STATE $b_{ij} = 1$
                \STATE $\overline{\beta^{\leq}_{ij}} = \overline{w_{j}^\leq}$
                \STATE $\overline{\beta^{\geq}_{ij}} = 0$
            \ELSE
                \STATE $b_{ij} = 0$
                \STATE $\overline{\beta^{\leq}_{ij}} = 0$
                \STATE $\overline{\beta^{\geq}_{ij}} = \overline{w_{j}^\geq}$
            \ENDIF
        \ENDFOR
    \ELSE
        \STATE $\overline{w_{j}^\leq} = 0$
        \STATE $\overline{w_{j}^\geq} = 0$
        \STATE $\forall i = 1,\dots,N$ $\overline{\beta^{\geq}_{ij}} = 0$
        \STATE $\forall i = 1,\dots,N$ $\overline{\beta^{\leq}_{ij}} = 0$
    \ENDIF
    \RETURN $\overline{\alpha}, \overline{\beta^{\leq}}, \overline{\beta^{\geq}}, \overline{b}, \overline{w_0}, \overline{w^\leq}, \overline{w^\geq}, \overline{t}$
\ENDFOR
\end{algorithmic}
\caption{A simple heuristic}\label{alg:heur}
\end{algorithm}

The first step is to solve the continuous relaxation of formulation \eqref{eq:linear_obj}-\eqref{eq:linear_last}, and obtain one of its optimal solutions, say $\overline{\alpha}, \overline{\beta^{\leq}}, \overline{\beta^{\geq}}, \overline{b}, \overline{w_0}, \overline{w^\leq}, \overline{w^\geq}, \overline{t}$. Then, we round $\overline{w_0}$ to the closest integer (Step \ref{algo1:round_w0}), and make $\overline{\alpha}$ binary, thus fixing to $1$ the coordinates of those predictors that will take part in the model (Step \ref{algo1:round_alphas}). To do so, consider $\overline{\alpha}$ in non-increasing order and round up to $1$ the first $k$ values, stopping when a $0$ is found. According to this, for each feature $j=1,\ldots,p$, we fix $\overline{w^\leq}_j, \overline{w^\geq}_j$ to the closest integer as long as $\overline{\alpha}_j$ takes value $1$ after the rounding of Step 3; otherwise, $\overline{w^\leq}_j= \overline{w^\geq}_j=0$. The value of $\overline{t}$ remains the same. Then, for each feature $j=1,\ldots,p$ and each individual $i=1,\ldots,N$, we assign $\overline{\beta_{ij}^{\leq}}, \overline{\beta_{ij}^{\geq}}$ the values of $\overline{w^\leq}_j,\overline{w^\geq}_j$ when $\overline{\alpha}_j = 1$ and $x_{ij}\leq \overline{t}_j$, $x_{ij} > \overline{t}_j$, respectively, and $0$, otherwise. At the same time, $b_{ij}$ is also fixed to $1$ and $0$, depending if $x_{ij}\leq\overline{t}_j$ or not.

\section{Computational experience}
\label{Computational experience}

The purpose of this section is to compare the performance of our formulation with state-of-the-art MINLO solvers against Algorithm \ref{alg:heur}. To this aim, we investigate the computational effort required to retrieve an optimal solution, as a function of $p$, $k$ and $N$. 

\subsection{Setup}
\label{Setup}


The Mixed-Integer Non-Linear Optimization formulation \eqref{eq:linear_obj}-\eqref{eq:linear_last} has been implemented using Pyomo optimization modeling language \cite{bynum2021pyomo} in Python
3.11 \cite{python}. We have tested BONMIN 1.8.9 solver using CBC 2.10.8 and IPOPT 3.12.13, with the algorithm \texttt{B-Hyb}; and Algorithm \ref{alg:heur} using IPOPT 3.12.13 solver as input. The time limit was set to 1000s in both cases. Our experiments have been conducted on a PC, with an 13th GenIntel$^\circledR$ Core$^{\rm TM}$ i7-1365U CPU 1.80GHz processor and 64 GB RAM. The operating system is 64 bits. The code can be found at https://github.com/mmolerous/Risk-Scores.

We have designed synthetic data sets with $N$ ranging from $200$ to $20000$ observations, and $p \in \left\lbrace 1,2,3,4,5 \right\rbrace$ continuous predictors, which were generated following continuous univariate distributions in the unit interval, i.e., $X_1,\ldots,X_p \sim U\left[0,1\right]$. Random integer coefficients $w_0, w_1^\leq, \ldots, w_p^\leq, w_1^\geq, \ldots, w_p^\geq $ in $\left[ -5,5\right]$ and random thresholds $t_1,\ldots,t_p$ in $\left[0,1\right]$ were assigned to construct the scoring system. Then, the labels $\left\lbrace y_i\right\rbrace_{i=1}^{N}$  were assigned according to their probability value, see Equation \eqref{eq:prob1}. If $P(y_i=1|\bm{x}_i)\geq 0.5$, we assign the label $1$, and $-1$ otherwise. 

Four metrics are assessed for each of the experiments found in Tables \ref{tab:varying_p}, \ref{tab:varying_k}, \ref{tab:varying_N} and \ref{tab:varying_N_p2}. TT, total time spent by the solver; ST, time spent by the solver to attain the last best feasible solution found; OBJ, training logistic loss corresponding to the solution obtained in ST; and AUC, training area under the curve corresponding to the solution obtained in ST.
$\ast$ will denote that the time limit has been reached. $-$ will denote that the solver has even not able to generate a feasible solution within the time limit. Note that these tables contain the actual OBJ and AUC values for the synthetic datasets, under the Reference label.


\subsection{Results with respect to p}
\label{Results with respect to p}

$N=200$ is tested with respect to $p \in \left\lbrace 1,2,3,4,5 \right\rbrace$. The problem was solved for parameter $k=p$. Table \ref{tab:varying_p} shows the results. We observe that, for $p=1, 2$, BONMIN solver is able to retrieve the reference solution, within $10$s and $621$s, respectively, obtaining the best AUC values. BONMIN solver starts to struggle for $p\geq 3$, where the time limit is reached, and also it has not retrieved the solution because its OBJ and AUC values are worse than the referenced ones. This highlights the computational burden of the considered problem using state-of-the-art MINLO solvers. Algorithm \ref{alg:heur} performs fast to find a first feasible solution, with AUC values above 0.70. In particular, contrary to BONMIN solver, Algorithm \ref{alg:heur} produces much better performance for $p\geq 3$, which makes it promising.

 \begin{table}[H]
    \centering
\begin{tabular}{r | c c | c c c c | c c c }
& \multicolumn{2}{c|}{Reference} & \multicolumn{4}{c|}{BONMIN} &  \multicolumn{3}{c}{Algorithm \ref{alg:heur}} \\ \hline
$p$ & OBJ & AUC & TT & ST & OBJ & AUC & TT/ST & OBJ & AUC \\ \hline
$1$ & 96.16& 1.00 & 10.91 & 10.91 & 1.34 & 1.00 & 2.35 & 81.65 & 0.92 \\ \hline
$2$ &24.80 & 1.00 & 621.46 & 621.46 & 0.47 &1.00&0.46 & 43.34 & 0.92 \\  \hline
$3$ & 45.54 &1.00 & * & 207.07& 644.00& 0.63&  0.68& 72.44& 0.81 \\ \hline
$4$ & 26.34 &1.00 & *& 10.96 &138.62 & 0.50& 0.80 &70.44 &  0.79 \\ \hline
$5$ &38.36 & 1.00 &* & 10.49 &138.62 &0.50 &1.01  &109.11 & 0.72 \\ 
\end{tabular}
    \caption{For $N=200$, as a function of $p$, and $k=p$: TT, total time spent by the solver; ST, time spent by the solver to attain the last best feasible solution found; OBJ, training logistic loss corresponding to the solution obtained in ST; and AUC, training area under the curve corresponding to the solution obtained in ST. $\ast$ denotes that the time limit has been reached.}
    \label{tab:varying_p}
\end{table}

\subsection{Results with respect to k}
\label{Results with respect to k}

$N=200$ and $p=5$ is tested with respect to $k \in \left\lbrace 1,2,3,4,5 \right\rbrace$. Table \ref{tab:varying_k} shows the results. We observe that, for BONMIN solver, the time limit is reached for all the configurations tested, also obtaining bad performance. This is not the case for Algorithm \ref{alg:heur}, which is always able to find a first feasible solution that is better in terms of both OBJ and AUC within 1s.

 \begin{table}[H]
    \centering
\begin{tabular}{r | c c | c c c c | c c c }
& \multicolumn{2}{c|}{Reference} & \multicolumn{4}{c|}{BONMIN} &  \multicolumn{3}{c}{Algorithm \ref{alg:heur}}  \\ \hline
$k$ & OBJ & AUC & TT & ST & OBJ & AUC & TT/ST & OBJ & AUC   \\ \hline
$1$ &38.36&1.00 & * & 10.69 & 138.62 & 0.50& 1.23&  124.76 & 0.79\\ \hline
$2$ &38.36 & 1.00 & * & 10.79& 138.62 & 0.50 & 1.23& 105.58 & 0.78\\ \hline
$3$ &38.36 & 1.00 & * & 10.50 & 138.62 & 0.50& 1.45 & 92.29 & 0.79 \\ \hline
$4$ & 38.36 & 1.00 & * & 13.03 & 138.62& 0.50& 1.58& 95.57 & 0.75 \\ \hline
$5$ &38.36 & 1.00 &* & 10.49 &138.62 &0.50 &1.01 &109.11& 0.72 \\ 
\end{tabular}
    \caption{For $N=200$ and $p=5$, as a function of $k$: TT, total time spent by the solver; ST, time spent by the solver to attain the last best feasible solution found; OBJ, training logistic loss corresponding to the solution obtained in ST; and AUC, training area under the curve corresponding to the solution obtained in ST. $\ast$ denotes that the time limit has been reached.
    }
    \label{tab:varying_k}
\end{table}

\subsection{Results with respect to N}
\label{Results with respect to N}

$p=1$ and $N \in \left\lbrace 200,1000,5000,10000,15000,20000 \right\rbrace$ were tested, solving for $k=p$. Table \ref{tab:varying_N} shows the results. For BONMIN solver, the time limit is reached for $N=20000$, although ST is always under $500$s, obtaining the best AUC values for all the configuration tested. 
Algorithm \ref{alg:heur} still produces feasible solutions quite fast up to $N=20000$, with AUC values around 0.90.

 \begin{table}[H]
    \centering
\begin{tabular}{r | c c | c c c c | c c c }
& \multicolumn{2}{c|}{Reference} & \multicolumn{4}{c|}{BONMIN} &  \multicolumn{3}{c}{Algorithm \ref{alg:heur}}\\ \hline
$N$ & OBJ & AUC & TT & ST & OBJ & AUC & TT/ST & OBJ & AUC\\ \hline
$200$ & 96.16& 1.00 & 10.91 & 10.91 & 1.34 & 1.00 &2.35 & 81.65 & 0.92 \\ \hline
$1000$ & 458.16 & 1.00 & 15.83 & 15.83& 6.71&  1.00 & 0.94 &  412.26& 0.91\\ \hline
$5000$ &2315.74 & 1.00 &57.92 &57.92 & 33.57& 1.00  &5.21  & 2050.30 & 0.91 \\ \hline
$10000$ & 4650.17 & 1.00 & 144.31 & 144.31 &67.15 & 1.00  &9.14 & 4095.61 & 0.91  \\ \hline
$15000$ & 6949.49&1.00 & 405.31&  405.31 &100.73 & 1.00 & 13.78& 6105.92 & 0.92 \\ \hline
$20000$ & 9276.00 & 1.00 & *  & 446.75  & 404.45& 1.00 & 18.75 &8145.23 & 0.92  \\ 

\end{tabular}
    \caption{For $p,k=1$, as a function of $N$: TT, total time spent by the solver; ST, time spent by the solver to attain the last best feasible solution found; OBJ, training logistic loss corresponding to the solution obtained in ST; and AUC, training area under the curve corresponding to the solution obtained in ST. $\ast$ denotes that the time limit has been reached.}
        \label{tab:varying_N}
\end{table}

A similar experiment was done for $p=2$ and $N \in \left\lbrace 200, 1000, 5000, 10000, 15000, 20000\right\rbrace$, solving for $k=p$. Table \ref{tab:varying_N_p2} shows the results. For BONMIN solver, the time limit is reached from $N=1000$ onwards. In particular, for $N\geq 5000$, it is not even able to produce a feasible solution. With respect to Algorithm \ref{alg:heur}, it is able to obtain a first feasible solution in less than $50$s for all the configurations tested, with AUC values above 0.90. 

 \begin{table}[H]
    \centering
\begin{tabular}{r | c c | c c c c | c c c }
& \multicolumn{2}{c|}{Reference} & \multicolumn{4}{c|}{BONMIN} &  \multicolumn{3}{c}{Algorithm \ref{alg:heur}} \\ \hline
$N$ & OBJ & AUC & TT & ST & OBJ & AUC & TT/ST & OBJ & AUC\\ \hline
$200$ &24.80 & 1.00 & 621.46 & 621.46 & 0.47 &1.00& 0.46 &43.34 & 0.92  \\ \hline
$1000$ & 118.00& 1.00 & * & 730.72 & 102.29 & 0.97 & 2.27 &  232.39 & 0.92  \\ \hline
$5000$ & 577.20 & 1.00 & * & - & - & - & 11.90& 1170.77 & 0.92  \\ \hline
$10000$ & 1163.40 & 1.00 & * & - & - & - & 21.65 & 2316.01& 0.91 \\ \hline
$15000$ & 1694.82 &  1.00& * & - & - & - &33.59  & 3468.51 & 0.91  \\ \hline
$20000$ & 9276.00 & 1.00 & * & - & - & - & 44.97& 4596.68 & 0.92 \\ 
\end{tabular}
    \caption{For $p=2$, as a function of $N$: TT, total time spent by the solver; ST, time spent by the solver to attain the last best feasible solution found; and ACC, training accuracy corresponding to the solution obtained in ST. $\ast$ denotes that the time limit has been reached. $-$ denotes that the solver has even not able to generate a feasible solution within the time limit.}
        \label{tab:varying_N_p2}
\end{table}

\section{Conclusions and perspectives}
\label{Conclusions and perspectives}

Risk Scores are one of the fundamental forms of Interpretable Machine Learning. In recent years, several approaches have been proposed to construct risk scores that handle binary data. In this paper, we give a Mixed-Integer Non-Linear Optimization formulation with linear constraints to build a risk score for continuous predictors. The model is tested in synthetic datasets.

Our findings illustrate how quickly the problem gets time-consuming for state-of-the-art MINLO solvers. To overcome this computational burden, a first simple matheuristic is presented, i.e., mathematical model-based heuristics. This matheuristic is shown to produce fast, and sometimes competitive, feasible solutions for all the configurations tested, which seems promising. Then, different steps are foreseen to follow. First, the output of this matheuristic could be used as initial solution of a more sophisticated one in order to improve prediction performance and scale up with the size of the dataset. Second, the definition of new constraints that may strengthen the formulation is also interesting. Last, one could linearize the objective function and focus on solving a Mixed Integer Linear Optimization (MILO) problem instead. Since there exist practically efficient solution algorithms for MILO problems, one could obtain relatively fast a solution for the MILO problem that will still be feasible for the MINLO problem since the constraints are linear.


\section*{Acknowledgments}

This publication was supported by the Chair ``Integrated Urban Mobility'', backed by L’X - \'Ecole Polytechnique and La Fondation de l’\'Ecole Polytechnique. The Partners of the Chair shall not under any circumstances accept any liability for the content of this publication, for which the author shall be solely liable. The support from research projects 
FQM-329 (Junta de Andalucía, Spain) and PID2022-137818OB-I00 (funded by MCIN/AEI/10.13039/501100011033, Spain, and FEDER/UE) are also gratefully acknowledged.

%
%
%
\bibliographystyle{splncs04}
\bibliography{OptimalRiskScoresForContinuousPredictors.bib}

\end{document}